\documentclass[11pt]{article}
\usepackage{amsfonts}
\usepackage{latexsym,amsmath}
\usepackage{amssymb,array}
\usepackage[authoryear,round]{natbib}
\usepackage{setspace}
\usepackage{enumitem}
\usepackage{xcolor}
\setlist{noitemsep}
\makeatletter \oddsidemargin 0in \evensidemargin 0in \textwidth 16cm
\usepackage{graphicx}
\providecommand{\keywords}[1] {
  \textit{Keywords:} #1}
\begin{document}
    \newcommand{\ov}{\overline}
    \newcommand{\om}{\omega}
    \newcommand{\ga}{\gamma}
    \newcommand{\cd}{\circledast}
    \newtheorem{thm}{Theorem}[section]
    \newtheorem{remark}{Remark}[section]
    \newtheorem{counterexample}{Counterexample}[section]
    \newtheorem{coro}{Corollary}[section]
    \newtheorem{propo}{Proposition}[section]
    \newtheorem{definition}{Definition}[section]
    \newtheorem{example}{Example}[section]
    \newtheorem{lem}{Lemma}[section]
    %\numberwithin{equation}{section}
    \date{}

\title{Dispersive and star ordering of sample extremes from dependent random variables following the proportional odds model}
\author{Arindam Panja$^1$\footnote
    {Corresponding author, e-mail: arindampnj@gmail.com}, Pradip Kundu$^2$
     and Biswabrata Pradhan$^1$\\
$^1$ SQC \& OR Unit, Indian Statistical Institute, Kolkata-700108, India\\
$^2$ Decision Science and Operations Management, Birla Global
University,\\ Bhubaneswar, Odisha-751003, India\\\footnote{Email
Address: kundu.maths@gmail.com (Pradip Kundu), bis@isical.ac.in
(Biswabrata Pradhan)}} \maketitle
\begin{abstract}
Dispersive order is a type of variability order for comparing the
variability in probability distributions. Star order compares the
skewness of probability distributions. This work considers
dispersive and star orders of extreme order statistics from
dependent random variables following the proportional odds (PO)
model. The joint distribution of the random variables is modeled
with Archimedean copula. Numerical examples are provided to
illustrate the findings.
\end{abstract}


\section{Introduction}
Suppose $X_{k:n}$, $k=1,2,\ldots,n$ denotes the $k$th order
statistic corresponding to random variables (r.v.'s)
$X_{1},X_{2},\ldots,X_{n}$. Order statistics play a crucial role in
statistical inference, reliability theory, life-testing, operations
research and economics. For example, in reliability theory, the
smallest and the largest order statistics $X_{1:n}$ and $X_{n:n}$,
respectively, represent the lifetimes of the series and the parallel
systems, where the corresponding r.v.'s represent the lifetimes of
$n$ components. Stochastic ordering has been widely used to compare
the magnitude and variability of extreme order statistics. However,
despite the importance and wide applications of the variability
orders (e.g. dispersive order and star order), there are less
research works in this direction as compared to the magnitude orders
(e.g., stochastic order, hazard rate order, reversed hazard rate
order, and likelihood ratio order).

Consider two random variables $X$ and $Y$ with cumulative distribution functions $F$ and $G$,denote $F^{-1}$ and $G^{-1}$ their respective right continuous
inverses, and of $\bar{F}$ and $\bar{G}$ their respective survival functions.
Then $X$ is said to be smaller than $Y$ in the
    \begin{enumerate}[label=(\roman*)]
        \item dispersive order (denoted as $X \leq_{disp} Y$)
        if $F^{-1}(v)-F^{-1}(u)\leq G^{-1}(v)-G^{-1}(u)$ for all $0\leq
        u\leq v\leq 1$. This is equivalently to $G^{-1}(F(x))-x$ is increasing
        in $x$. When $X$ and $Y$ have probability density functions(pdfs)
        $f$ and $g$, respectively, then $X \leq_{disp} Y$ if, and only if,
        $g(G^{-1}(u))\leq f(F^{-1}(u))$ for all $u\in (0,1)$ \citep{shaked};
        \item star order (denoted by $X \leq_{\star} Y$) if $G^{-1}(F(x))/x$ is increasing in $x\in \mathbb{R}_+$ \citep{shaked} .
    \end{enumerate}

Dispersive order is one kind of variability order for comparing
variability in probability distributions \citep{joen,skoch,shaked}.
Star order have been introduced in the literature to compare the
skewness of probability distributions. The star order is also called
more IFRA (increasing failure rate in average) order. If one r.v. is
smaller than another in terms of star order, then this can be
interpreted as the former r.v. ages faster than the later in the
sense of the star ordering. For more discussion and applications,
see \cite{barl}, \cite{skoch} and \cite{zhang}. Skewed distributions
often serve as reasonable models for system lifetimes, auction
theory, insurance claim amounts, financial returns etc. and thus it
is of interest to compare skewness of probability distributions
\citep{wu}. Recently, there have been a number of works on
dispersive and star ordering of extreme order statistics of random
samples from different family of distributions
\citep{ding,fang1,fang2,koch2,koch3,lif,nadeb,zhang,zhang1}. There
are some research works on sample spacings also, like \cite{xu}
established dispersive and star ordering for sample spacing from
heterogeneous exponential.
 distributions

The proportional odds (PO) model \citep{benn,kirmani} is a very
important model in reliability theory and survival analysis. Let $X$
and $Y$ be two r.v.'s with cdfs $F$, $G$, and survival functions
$\bar{F}$, $\bar{G}$, respectively. If the r.v. $X$ denote a
survival time, then the odds function $\theta_{X}(t)$ defined by
$\theta_{X}(t)=\bar{F}(t)/F(t)$ represents the odds on surviving
beyond time $t$. The r.v.'s $X$ and $Y$ are said to satisfy PO model
if $\theta_{Y}(t)=\alpha~\theta_{X}(t)$ for all admissible $t$,
where $\alpha$ is a proportionality constant known as proportional
odds ratio. Then the survival functions of $X$ and $Y$ are related
as
\begin{equation}\label{poalt}\bar{G}(t)=\frac{\alpha
\bar{F}(t)}{1-\bar{\alpha}\bar{F}(t)},\end{equation} where
$\bar{\alpha}=1-\alpha$. From this relation, it can be observed that
the ratio of hazard rate functions becomes
$1/(1-\bar{\alpha}\bar{F}(t))$, so that the hazard ratio is
increasing (decreasing) for $\alpha>1$ ($\alpha<1$) and it converges
to unity as $t$ tends to $\infty$. This is in contrast to the
proportional hazard rate (PHR) model where this ratio remains
constant with time. The convergence property of hazard functions
makes the PO model reasonable in many practical applications as
discussed by \cite{benn}, \cite{coll}, \cite{kirmani}, \cite{lu} and
\cite{ross}. Also, the model (\ref{poalt}) with $0<\alpha <\infty$
provides a method of generating more flexible new family of
distributions known as Marshall-Olkin family of distributions or
Marshall-Olkin extended distributions \citep{cord1,marsh1}, from an
existing family of distributions. Extended Weibull distributions,
extended linear failure-rate distributions and extended generalized
exponential distributions are few examples those have been widely
studied in the literature. Thus, model (\ref{poalt}) has
implications both in terms of the PO model and in extending any
existing family of distributions to add flexibility in modeling.
This makes the PO model worth investigating.
%We will say that the r.v. $Y$ is following the PO model
%with baseline survival function $\bar{F}(\cdot)$ and parameter
%(proportionality constant) $\alpha$.

\par Let {\bf X}  = ($X_{1},X_{2},\ldots,X_{n}$) have
joint distribution function $F$ and joint survival function
$\bar{F}$. The marginal distribution function and survival function
of $X_i$ are  $F_i$ and $\bar{F}_i$, respectively, $i=1,2,\ldots,n$.
If there exist $C$, $\bar{C}:[0,1]^n\mapsto [0,1]$ such that
$F(x_1,\ldots,x_n)=C(F_1(x_1),\ldots,F_n(x_n))$ and
$\bar{F}(x_1,\ldots,x_n)=\bar{C}(\bar{F}_1(x_1),\ldots,\bar{F}_n(x_n))$
for all $x_i$, $i\in I_n$, then $C$ and $\bar{C}$ are called the
copula and survival copula respectively. If
$\varphi:[0,+\infty)\mapsto [0,1]$ with $\varphi(0)=1$ and
$\lim_{t\rightarrow +\infty} \varphi(t)=0$, then
$C(u_1,\ldots,u_n)=\varphi(\varphi^{-1}(u_1)+\ldots+\varphi^{-1}(u_n))=\varphi(\sum_{i=1}^n
\phi (u_i))$ for all $u_i\in(0,1]$, $i\in I_n$ is called an
Archimedean copula with generator $\varphi$ provided
$(-1)^{k}\varphi^{(k)}(t)\geq 0$, $k=0,1,\ldots,n-2$ and
$(-1)^{n-2}\varphi^{(n-2)}(t)$ is decreasing and convex for all
$t\geq 0$. Here $\phi=\varphi^{-1}$ is the right continuous inverse
of $\varphi$ so that $\phi(u)=\varphi^{-1}(u)=\sup\{t\in
\mathbb{R}:\varphi(t)>u\}$. In case of dependent samples, Li and
Fang \cite{lif} derived the dispersive order between maximums of two
PHR samples having a common Archimedean copula. For samples
following scale model, \cite{lil} obtained the dispersive and the
star orders between minimums of one heterogeneous and one
homogeneous samples sharing a common Archimedean copula.
\cite{fang1} investigated the dispersive order and the star order of
extreme order statistics for the samples following PHR model with
Archimedean survival copulas. \cite{fang2} obtained the dispersive
order between minimums of two scale proportional hazards samples
with a common Archimedean survival copula. With resilience-scaled
components, \cite{zhang1} derived the dispersive and the star order
between parallel systems, one consisting dependent heterogeneous
components and another consisting homogeneous components sharing a
common Archimedean survival copula.

In case of PO model, some authors, e.g. \cite{kundu1},
\cite{kundu2}, \cite{panja}, \cite{nanda} have investigated
stochastic comparison of this family of distributions and sample
extreme in the sense of magnitude orders. To the best of our
knowledge, there is no related study on the variability of extreme
order statistics arising from independent or dependent r.v.'s
following the PO model. Motivated by this, in this paper, we develop
the dispersive and the star ordering for comparing the minimums and
the maximums of dependent samples following the PO model. The
organization of the rest of the paper is as follows. In Section 2,
we consider comparisons of minimum order statistics from dependent
samples following the PO model in terms of the dispersive order and
the star order. Section 3 investigates the comparison of maximum
order statistics in terms of dispersive and star orders. In Section
4, some examples are provided to illustrate the main results of the
paper. In Section 5, we make concluding remarks.

\section{The dispersive ordering and the star ordering of minimums of dependent samples following the PO model}

In this section we consider the dispersive ordering of minimums of dependent random variables. We compare stochastically the minimums of two dependent
samples, one formed from heterogeneous r.v.'s and another from
homogeneous r.v.'s. Let $X=(X_1,X_2,...,X_n)$ be a set of
dependent r.v.'s coupled with Archimedean survival copula with
generator $\varphi$ and following the PO model with baseline
survival function $\bar{F}$, denoted as $X\sim
PO(\bar{F},\boldsymbol\alpha,\varphi)$, where
$\boldsymbol\alpha=(\alpha_1,\alpha_2,...,\alpha_n)\in
\mathbb{R}^n_{+}$ is the proportional odds ratio vector. That is,
odds function of each r.v. $X_i$ is proportional to an odds function
(baseline odds) of a r.v. having distribution function $F$, with
proportionality constant $\alpha_i$. %Let $Y=(Y_1,Y_1,...,Y_1)$ and
%$Y\sim PO(\bar{F},\alpha\textbf{1},\varphi)$, $\alpha>0$, i.e. odds
%function of $Y_1$ is proportional to the baseline odds, with
%proportionality constant $\alpha$. The reliability functions
%(survival functions) of $X_i$ and $Y_1$ are
%$$\bar{F}_{X_i}(x)=\frac{\alpha_i
%\bar{F}(x)}{1-\bar{\alpha}_i\bar{F}(x)} ~\text{and}~
%\bar{F}_{Y_1}(x)=\frac{\alpha
%\bar{F}(x)}{1-\bar{\alpha}\bar{F}(x)}$$ where
%$\bar{\alpha}_i=1-\alpha_i$, $i=1,\ldots,n$ and
%$\bar{\alpha}=1-\alpha$.
We have the survival functions of $X_{1:n}$ as
\begin{equation}\label{eqspo}\bar{F}_{X_{1:n}}(x)=\varphi\left(\sum_{i=1}^n
\phi\left(\bar{F}_{X_i}(x)\right)\right),\end{equation} where
$\bar{F}_{X_i}(x)=\frac{\alpha_i\bar{F}(x)}{1-\bar{\alpha}_i\bar{F}(x)}$,
$\phi(u)=\varphi^{-1}(u)$, $u\in(0,1]$.

The following theorem consider the comparison of minimums of two samples, one from $n$ dependent
heterogeneous r.v.'s following the PO model and another from $n$ dependent homogeneous r.v.'s following the PO model, in terms of dispersive order. The result holds for the decreasing failure rate (DFR) baseline distribution $F$. The distribution function $F$ is said to be DFR  if the corresponding hazard rate $r(\cdot)$ is decreasing and increasing failure rate (IFR) distribution if $r(\cdot)$ is increasing.

\begin{thm}\label{thsdis} Suppose $X\sim PO(\bar{F},\boldsymbol\alpha,\varphi)$ and $Y\sim
PO(\bar{F},\alpha\textbf{1},\varphi)$. Then $X_{1:n}\leq_{disp}Y_{1:n}$ if the baseline distribution $F$ is DFR,
$\varphi$ is log-convex, $\frac{\varphi}{\varphi'}$ is concave and $\alpha\geq\frac{1}{n}\sum_{i=1}^n\alpha_{i}$, for $0\leq\alpha \leq 1$.
\end{thm}
\textbf{Proof:} We have the distribution functions of $X_{1:n}$ and
$Y_{1:n}$ as $F_{1}(x)=1-\varphi\left(\sum_{i=1}^n
\phi\left(\bar{F}_{X_i}(x)\right)\right)$ and
$G_1(x)=1-\varphi\left(n \phi\left(\bar{F}_{Y_1}(x)\right)\right)$,
respectively,  where
$\bar{F}_{X_i}(x)=\frac{\alpha_i\bar{F}(x)}{1-\bar{\alpha}_i\bar{F}(x)}$
and $\bar{F}_{Y_1}(x)=\frac{ \alpha
\bar{F}(x)}{1-\bar{\alpha}\bar{F}(x)}$, $x\in \mathbb{R}$. The
respective pdfs of $X_{1:n}$ and $Y_{1:n}$
are given by
\begin{equation}\label{eqf1s} f_1(x)=\varphi'\left(\sum_{i=1}^n
\phi\left(\bar{F}_{X_i}(x)\right)\right)\sum_{i=1}^n
\frac{\varphi\left(\phi\left(\bar{F}_{X_i}(x)\right)\right)}{\varphi'\left(\phi\left(\bar{F}_{X_i}(x)\right)\right)}
\frac{r(x)}{1-\bar{\alpha}_i \bar{F}(x)},
\end{equation}
and
\begin{equation*} g_1(x)=n \varphi'\left(n
    \phi\left(\bar{F}_{Y_1}(x)\right)\right)\cdot
    \frac{r(x)}{1-\bar{\alpha} \bar{F}(x)}\cdot
    \frac{\varphi\left(\phi\left(\bar{F}_{Y_1}(x)\right)\right)}{\varphi'\left(\phi\left(\bar{F}_{Y_1}(x)\right)\right)},
\end{equation*}
We have
\begin{equation*} G_1^{-1}(x)=\bar{F}^{-1}\left(\frac{
\varphi\left(\frac{1}{n}\phi(1-x)\right)}{\alpha+\bar{\alpha}\varphi\left(\frac{1}{n}\phi(1-x)\right)}\right).
\end{equation*}
So
\begin{equation}\label{eqGFs} G_1^{-1}(F_{1}(x))=\bar{F}^{-1}\left(\frac{
\varphi\left(\frac{1}{n}\sum_{i=1}^n
\phi\left(\bar{F}_{X_i}(x)\right)\right)}{\alpha+\bar{\alpha}\varphi\left(\frac{1}{n}\sum_{i=1}^n
\phi\left(\bar{F}_{X_i}(x)\right)\right)}\right)
=\bar{F}^{-1}(\gamma(x)),
\end{equation} where $\gamma(x)=\frac{
\varphi\left(\frac{1}{n}\sum_{i=1}^n
\phi\left(\bar{F}_{X_i}(x)\right)\right)}{\alpha+\bar{\alpha}\varphi\left(\frac{1}{n}\sum_{i=1}^n
\phi\left(\bar{F}_{X_i}(x)\right)\right)}$.\\
Now
\begin{eqnarray}\nonumber g_{1}(G_1^{-1}(F_{1}(x)))&=&n \varphi'\left(n
\phi\left(\frac{ \alpha
\gamma(x)}{1-\bar{\alpha}\gamma(x)}\right)\right)\cdot
\frac{r\left(\bar{F}^{-1}(\gamma(x))\right)}{1-\bar{\alpha}
\gamma(x)}\cdot \frac{\varphi\left(\phi\left(\frac{ \alpha
\gamma(x)}{1-\bar{\alpha}\gamma(x)}\right)\right)}{\varphi'\left(\phi\left(\frac{
\alpha\gamma(x)}{1-\bar{\alpha}\gamma(x)}\right)\right)}\\
\nonumber &=&n \varphi'\left(\sum_{i=1}^n
\phi\left(\bar{F}_{X_i}(x)\right)\right)\cdot \left(\alpha
+\bar{\alpha}\varphi\left(\frac{1}{n}\sum_{i=1}^n
\phi\left(\bar{F}_{X_i}(x)\right)\right)\right) \nonumber \\
& & \times \frac{\varphi\left(\frac{1}{n}\sum_{i=1}^n
\phi\left(\bar{F}_{X_i}(x)\right)\right)}{\varphi'\left(\frac{1}{n}\sum_{i=1}^n
\phi\left(\bar{F}_{X_i}(x)\right)\right)}\cdot
\frac{r\left(\bar{F}^{-1}(\gamma(x))\right)}{\alpha}.\label{eqgGs}  \end{eqnarray}
Note that $\bar{F}_{X_i}(x)$ is increasing and concave in $\alpha_i$
and $1/(1-\bar{\alpha}_i \bar{F}(x))$ is decreasing and convex in
$\alpha_i$. Also it can be seen that
$\phi\left(\bar{F}_{X_i}(x)\right)$ is decreasing and convex in
$\alpha_i$ if $\varphi$ is log-convex. Now denote
$\frac{1}{n}\sum_{i=1}^n\alpha_{i}=\boldsymbol\alpha^{avg}$ and
$\eta(\alpha_i)=\phi\left(\bar{F}_{X_i}(x)\right)$. Then for
$\alpha\geq\frac{1}{n}\sum_{i=1}^n\alpha_{i}=\boldsymbol\alpha^{avg}$,
from the convexity and decreasing property of
$\eta(\alpha_i)=\phi\left(\bar{F}_{X_i}(x)\right)$ with respect to
$\alpha_i$, we have
$\frac{1}{n}\sum_{i=1}^n\eta(\alpha_{i})\geq\eta(\boldsymbol\alpha^{avg})\geq\eta(\alpha)$,
which gives
\begin{eqnarray}
\label{conphi}\frac{1}{n}\sum_{i=1}^n
\phi\left(\bar{F}_{X_i}(x)\right)&\geq&
\phi\left(\bar{F}_{Y_1}(x)\right)  \\
\implies
\frac{\alpha}{\bar{\alpha}} +\varphi\left(\frac{1}{n}\sum_{i=1}^n
\phi\left(\bar{F}_{X_i}(x)\right)\right) &\leq&
\frac{\alpha}{\bar{\alpha}}+ \bar{F}_{Y_1}(x) \nonumber\\
\implies 1-\frac{\frac{\alpha}{\bar{\alpha}}}{\frac{\alpha}{\bar{\alpha}}+\varphi\left(\frac{1}{n}\sum_{i=1}^n
\phi\left(\bar{F}_{X_i}(x)\right)\right)} &\leq&
1-\frac{\frac{\alpha}{\bar{\alpha}}}{\frac{\alpha}{\bar{\alpha}}+\bar{F}_{Y_1}(x)}
\nonumber \\
\implies \frac{\varphi\left(\frac{1}{n}\sum_{i=1}^n
\phi\left(\bar{F}_{X_i}(x)\right)\right)}{\alpha+\bar{\alpha}\varphi\left(\frac{1}{n}\sum_{i=1}^n
\phi\left(\bar{F}_{X_i}(x)\right)\right)} &\leq&
\frac{\bar{F}_{Y_1}(x)}{\alpha+\bar{\alpha}\bar{F}_{Y_1}(x)}. \nonumber
\end{eqnarray}
This implies $\gamma(x) \leq \bar{F}(x).$ As a result we have $\bar{F}^{-1}(\gamma(x))\geq x$. Now if $r(\cdot)$ is
decreasing then
\begin{equation}\label{eqth11s}r(\bar{F}^{-1}(\gamma(x)))\leq r(x).
\end{equation}
Now from (\ref{conphi}), we have
\begin{eqnarray}\nonumber \alpha+\bar{\alpha}\varphi\left(\frac{1}{n}\sum_{i=1}^n
\phi\left(\bar{F}_{X_i}(x)\right)\right) &\leq&
\alpha+\bar{\alpha}\bar{F}_{Y_1}(x)\\
\nonumber &=& \frac{\alpha}{1-\bar{\alpha} \bar{F}(x)}\\
\label{eqth13s}&\leq& \alpha
\frac{1}{n}\sum_{i=1}^{n}\frac{1}{1-\bar{\alpha_i} \bar{F}(x)},
\end{eqnarray}
where the last inequality follows from the fact that
$\frac{1}{1-\bar{\alpha_i} \bar{F}(x)}$ is decreasing and convex in
$\alpha_i$.\\ If $\frac{\varphi}{\varphi'}$ is concave, then we have
\begin{equation}-\frac{\varphi\left(\frac{1}{n}\sum_{i=1}^n
\phi\left(\bar{F}_{X_i}(x)\right)\right)}{\varphi'\left(\frac{1}{n}\sum_{i=1}^n
\phi\left(\bar{F}_{X_i}(x)\right)\right)}\leq
-\frac{1}{n}\sum_{i=1}^n
\frac{\varphi\left(\phi\left(\bar{F}_{X_i}(x)\right)\right)}{\varphi'\left(\phi\left(\bar{F}_{X_i}(x)\right)\right)}.\end{equation}
Thus we have
\begin{equation}\label{alpha+}\left(\alpha+\bar{\alpha}\varphi\left(\frac{1}{n}\sum_{i=1}^n
\phi\left(\bar{F}_{X_i}(x)\right)\right)\right)\left(-\frac{\varphi\left(\frac{1}{n}\sum_{i=1}^n
\phi\left(\bar{F}_{X_i}(x)\right)\right)}{\varphi'\left(\frac{1}{n}\sum_{i=1}^n
\phi\left(\bar{F}_{X_i}(x)\right)\right)}\right)\leq \alpha
\frac{1}{n}\sum_{i=1}^n\left(-
\frac{\varphi\left(\phi\left(\bar{F}_{X_i}(x)\right)\right)}{\varphi'\left(\phi\left(\bar{F}_{X_i}(x)\right)\right)}\right)\frac{1}{n}\sum_{i=1}^{n}\frac{1}{1-\bar{\alpha_i}
\bar{F}(x)}
\end{equation}
If $\varphi$ is log-convex, then $-\frac{\varphi(x)}{\varphi'(x)}$
is increasing in $x$, so that
$-\frac{\varphi\left(\phi\left(\bar{F}_{X_i}(x)\right)\right)}{\varphi'\left(\phi\left(\bar{F}_{X_i}(x)\right)\right)}$
is decreasing in $\alpha_i$. So by Chebyshev's inequality we have
\begin{equation}\label{cheby}\frac{1}{n}\sum_{i=1}^n\left(-
\frac{\varphi\left(\phi\left(\bar{F}_{X_i}(x)\right)\right)}{\varphi'\left(\phi\left(\bar{F}_{X_i}(x)\right)\right)}\right)\cdot\frac{1}{n}\sum_{i=1}^{n}\frac{1}{1-\bar{\alpha_i}
\bar{F}(x)}\leq \frac{1}{n}\sum_{i=1}^n \left(-
\frac{\varphi\left(\phi\left(\bar{F}_{X_i}(x)\right)\right)}{\varphi'\left(\phi\left(\bar{F}_{X_i}(x)\right)\right)}\right)\frac{1}{1-\bar{\alpha_i}
\bar{F}(x)}
\end{equation}
From (\ref{eqth11s}), (\ref{alpha+}), (\ref{cheby}) and the fact
that the common factor $\varphi'\left(\sum_{i=1}^n
\phi\left(\bar{F}_{X_i}(x)\right)\right)$ in (\ref{eqf1s}) and
(\ref{eqgGs}) is negative, we have $g_{1}(G_1^{-1}(F_{1}(x)))\leq
f_{1}(x)$ for all $x\in \mathbb{R}$. Hence the theorem follows.$\hfill\Box$\\

It may be interest to know whether as in case of Theorem \ref{thsdis} we can
establish dispersive ordering for $\alpha\geq1$ when the baseline
distribution is IFR or DFR. The following counterexample shows that
with these conditions, we cannot establish dispersive ordering even
in case of samples from independent r.v.'s.

\begin{counterexample}\label{countex1}
Consider the minimums of two samples, one having three independent and
heterogeneous r.v.'s, and another having three independent and
homogeneous r.v.'s with respective distribution functions
$F_1(x)=1-\prod_{i=1}^{3}\left(\frac{\alpha_i\bar{F}(x)}{1-\bar{\alpha}_i\bar{F}(x)}\right)$
and $G_{1}(x)=1-\left(\frac{ \alpha
\bar{F}(x)}{1-\bar{\alpha}\bar{F}(x)}\right)^3$, where $\alpha_1=7$,
$\alpha_2=25$, $\alpha_3=100$,
$\alpha=(\alpha_1+\alpha_2+\alpha_3)/3=44$, and
$\bar{F}(x)=e^{-(9x)^{0.9}}$, so that the baseline distribution is
DFR. We obtain
\begin{equation*}g_{1}(G_1^{-1}(F_{1}(x)))=\frac{1}{\alpha}~3 \left(\prod_{i=1}^{3}
\bar{F}_{X_i}(x)\right)
\left(\alpha+\bar{\alpha}\left(\prod_{i=1}^{3}
\bar{F}_{X_i}(x)\right)^{1/3}\right)r(\bar{F}^{-1}(\gamma(x))),
\end{equation*}
where $\gamma(x)=\frac{ \left(\prod_{i=1}^{3}
\bar{F}_{X_i}(x)\right)^{1/3}}{\alpha+\bar{\alpha}\left(\prod_{i=1}^{3}
\bar{F}_{X_i}(x)\right)^{1/3}}$,\\
\noindent and
\begin{equation*}f_1(x)=\left(\prod_{i=1}^{3} \bar{F}_{X_i}(x)\right) r(x)
\left(\sum_{i=1}^3 \frac{1}{1-\bar{\alpha_i}\bar{F}(x)}\right).
\end{equation*}
We plot $g_{1}(G_1^{-1}(F_{1}(x)))-f_{1}(x)$ by substituting
$x=t/(1-t)$, so that for $x\in[0,\infty)$, we have $t\in[0,1)$. The plot is shown in Figure \ref{Figsdisp}(a) we observe from the plot that
$g_{1}(G_1^{-1}(F_{1}(x)))\nleq f_{1}(x)$ and also
$g_{1}(G_1^{-1}(F_{1}(x)))\ngeq f_{1}(x)$.\\
Next we take $\alpha_1=0.78$, $\alpha_2=0.97$, $\alpha_3=67$,
$\alpha=(\alpha_1+\alpha_2+\alpha_3)/3=22.9167$, and
$\bar{F}(x)=e^{-x^3}$, so that the baseline distribution is IFR. Figure \ref{Figsdisp}(b) illustrates the
plot of $g_{1}(G_1^{-1}(F_{1}(x)))-f_{1}(x)$ by substituting
$x=t/(1-t)$, so that for $x\in[0,\infty)$, we have $t\in[0,1)$. From
Figure \ref{Figsdisp}(b) we observe that
$g_{1}(G_1^{-1}(F_{1}(x)))-f_{1}(x)\nleq 0$ and also
$g_{1}(G_1^{-1}(F_{1}(x)))-f_{1}(x)\ngeq 0$.
\begin{figure}
\begin{center}
\includegraphics[width=12cm]{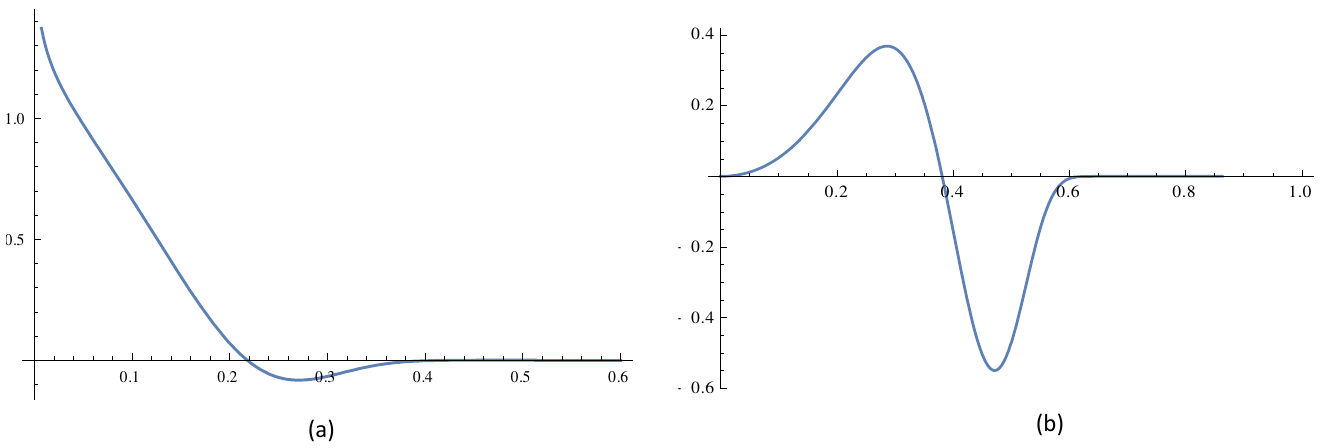}
\caption{Plot of $g_{1}(G_1^{-1}(F_{1}(x)))-f_{1}(x)$ for
$x=t/(1-t)$, $t\in[0,1]$ when baseline distribution is (a) DFR and
(b) IFR.}\label{Figsdisp}
\end{center}
\end{figure}
\end{counterexample}
The following theorem compare the minimums of two samples, both from
$n$ dependent homogeneous r.v.'s following the PO model and with
different Archimedean copulas.
\begin{thm}\label{thsdis2} Suppose $X\sim PO(\bar{F},\alpha\textbf{1},\varphi_1)$ and $Y\sim
PO(\bar{F},\alpha\textbf{1},\varphi_2)$. Then
$X_{1:n}\leq_{disp}Y_{1:n}$ if the baseline distribution is DFR,
$\varphi_2(\phi_2(w)/n)/\varphi_1(\phi_1(w)/n)$ is increasing in $w$
and $0\leq\alpha \leq 1$.
\end{thm}
\textbf{Proof:} The distribution functions of $X_{1:n}$ and
$Y_{1:n}$ are given by $G_{1}(x)=1-\varphi_1\left(n
\phi_1\left(\bar{F}_{X_1}(x)\right)\right)$, and
$G_2(x)=1-\varphi_2\left(n
\phi_2\left(\bar{F}_{X_1}(x)\right)\right)$, respectively, where
$\bar{F}_{X_1}(x)=\frac{ \alpha
\bar{F}(x)}{1-\bar{\alpha}\bar{F}(x)}$, $x\in \mathbb{R}$. The
respective pds are given by
\begin{equation}\label{g1dis2} g_1(x)=n \varphi_1'\left(n
\phi_1\left(\bar{F}_{X_1}(x)\right)\right)\frac{\varphi_1\left(\phi_1\left(\bar{F}_{X_1}(x)\right)\right)}{\varphi_1'\left(\phi_1\left(\bar{F}_{X_1}(x)\right)\right)}\cdot
\frac{r(x)}{1-\bar{\alpha} \bar{F}(x)},
\end{equation}
\begin{equation*} g_2(x)=n \varphi_2'\left(n
\phi_2\left(\bar{F}_{X_1}(x)\right)\right)\frac{\varphi_2\left(\phi_2\left(\bar{F}_{X_1}(x)\right)\right)}{\varphi_2'\left(\phi_2\left(\bar{F}_{X_1}(x)\right)\right)}\cdot
\frac{r(x)}{1-\bar{\alpha} \bar{F}(x)}.
\end{equation*}
We have
\begin{equation}\label{eqG2G1s} G_2^{-1}(G_{1}(x))=\bar{F}^{-1}\left(\frac{
\varphi_2\left(\frac{1}{n}\phi_2\left(\varphi_1\left(n
\phi_1\left(\bar{F}_{X_1}(x)\right)\right)\right)\right)}{\alpha+\bar{\alpha}\varphi_2\left(\frac{1}{n}\phi_2\left(\varphi_1\left(n
\phi_1\left(\bar{F}_{X_1}(x)\right)\right)\right)\right)}\right)
=\bar{F}^{-1}(\eta(x)),
\end{equation}
where $\eta(x) = \frac{
    \varphi_2\left(\frac{1}{n}\phi_2\left(\varphi_1\left(n
    \phi_1\left(\bar{F}_{X_1}(x)\right)\right)\right)\right)}{\alpha+\bar{\alpha}\varphi_2\left(\frac{1}{n}\phi_2\left(\varphi_1\left(n
    \phi_1\left(\bar{F}_{X_1}(x)\right)\right)\right)\right)}$.
\begin{eqnarray}\nonumber g_{2}(G_2^{-1}(G_{1}(x)))&=&n \varphi_2'\left(\phi_2\left(\varphi_1\left(n
\phi_1\left(\bar{F}_{X_1}(x)\right)\right)\right)\right)\frac{\varphi_2\left(\frac{1}{n}\phi_2\left(\varphi_1\left(n
\phi_1\left(\bar{F}_{X_1}(x)\right)\right)\right)\right)}{\varphi_2'\left(\frac{1}{n}\phi_2\left(\varphi_1\left(n
\phi_1\left(\bar{F}_{X_1}(x)\right)\right)\right)\right)}\frac{r\left(\bar{F}^{-1}(\eta(x))\right)}{\alpha}\\
\label{g2ginvs}&& \times
\left(\alpha+\bar{\alpha}\varphi_2\left(\frac{1}{n}\phi_2\left(\varphi_1\left(n
\phi_1\left(\bar{F}_{X_1}(x)\right)\right)\right)\right)\right).\end{eqnarray}
From Lemma 3.9 of \cite{fang1}, for increasing
$\varphi_2(\phi_2(w)/n)/\varphi_1(\phi_1(w)/n)$ we have
$\varphi_2\left(n \phi_2\left(\bar{F}_{X_1}(x)\right)\right)\geq
\varphi_1\left(n \phi_1\left(\bar{F}_{X_1}(x)\right)\right)$, which
implies $\bar{F}_{X_1}(x) \geq
\varphi_2\left(\frac{1}{n}\phi_2\left(\varphi_1\left(n
\phi_1\left(\bar{F}_{X_1}(x)\right)\right)\right)\right)$. Again
from this we get $\bar{F}(x)\geq \eta(x)$ which implies
$\bar{F}^{-1}(\eta(x))\geq x$. Thus if $r(\cdot)$ is decreasing then
\begin{equation}\label{eqrinvs}r(\bar{F}^{-1}(\eta(x)))\leq r(x).
\end{equation} Also for $\bar{\alpha}\geq0$, $\varphi_2\left(\frac{1}{n}\phi_2\left(\varphi_1\left(n
\phi_1\left(\bar{F}_{X_1}(x)\right)\right)\right)\right)\leq
\bar{F}_{X_1}(x)$ implies
\begin{equation}\label{eqalp+alpb} \alpha+\bar{\alpha}\varphi_2\left(\frac{1}{n}\phi_2\left(\varphi_1\left(n
\phi_1\left(\bar{F}_{X_1}(x)\right)\right)\right)\right)\leq
\frac{\alpha}{1-\bar{\alpha} \bar{F}(x)}.
\end{equation}
Again from Lemma 3.9 of
\cite{fang1} by substituting
$w=\varphi_1\left(n \phi_1\left(\bar{F}_{X_1}(x)\right)\right)$ in
increasing $\varphi_1(\phi_1(w)/n)/\varphi_2(\phi_2(w)/n)$, we get
\begin{equation}\label{remain}\frac{\varphi_2'\left(\phi_2\left(\varphi_1\left(n
\phi_1\left(\bar{F}_{X_1}(x)\right)\right)\right)\right)
\varphi_2\left(\frac{1}{n}\phi_2\left(\varphi_1\left(n
\phi_1\left(\bar{F}_{X_1}(x)\right)\right)\right)\right)}{\varphi_2'\left(\frac{1}{n}\phi_2\left(\varphi_1\left(n
\phi_1\left(\bar{F}_{X_1}(x)\right)\right)\right)\right)}\leq
\frac{\varphi_1'\left(n\phi_1\left(\bar{F}_{X_1}(x)\right)\right)
\varphi_1\left(\phi_1\left(\bar{F}_{X_1}(x)\right)\right)}{\varphi_1'\left(\phi_1\left(\bar{F}_{X_1}(x)\right)\right)}.\end{equation}
Now using (\ref{eqrinvs}), (\ref{eqalp+alpb}) and (\ref{remain}), we
have $g_{2}(G_2^{-1}(G_{1}(x)))\leq g_1(x)$ for all $x\in
\mathbb{R}$. This completes the proof.$\hfill\Box$

\begin{remark}\normalfont It is to be noted that Archimedean copula with generators
\begin{itemize}
\item [(i)] $\varphi_1(x)=\left(1+x^{1/\theta_1}\right)^{-\theta_1}$ and $\varphi_2(x)=\frac{1}{\left(x^{1/\theta_2}+1\right)}$  where $\theta_1\in(1,\infty), \theta_2 \in (1,\infty)$ and $1<\theta_1<\theta_2<\infty$,
\item[(ii)] $\varphi_1(x)= \left(\theta_1 x+1\right)^{-1/\theta_1}$ and $\varphi_2(x)=\left[\theta_2 x+1\right]^{-1/\theta_2},$ for
$0<\theta_1<\theta_2<\infty$,
\end{itemize}
satisfy the condition that
$\varphi_2(\phi_2(w)/n)/\varphi_1(\phi_1(w)/n)$ is increasing in $w$
for all $n\in Z$.
\end{remark}

The following corollary follows from Theorems \ref{thsdis} and
\ref{thsdis2}. This corollary compare the minimums of two samples,
one from from $n$ dependent heterogeneous r.v.'s following the PO
model and another from $n$ dependent homogeneous r.v.'s following
the PO model and with different Archimedean copulas.
\begin{coro}Suppose $X\sim PO(\bar{F},\boldsymbol\alpha,\varphi_1)$ and $Y\sim
PO(\bar{F},\alpha\textbf{1},\varphi_2)$. Then for
$\alpha\geq\frac{1}{n}\sum_{i=1}^n\alpha_{i}$,
$X_{1:n}\leq_{disp}Y_{1:n}$ if the baseline distribution is DFR,
$\varphi_1$ is log-convex, $\frac{\varphi_1}{\varphi_1'}$ is
concave, $\varphi_2(\phi_2(w)/n)/\varphi_1(\phi_1(w)/n)$ is
increasing in $w$, and $0\leq\alpha \leq
1$.\end{coro}\textbf{Proof:} Let $Z\sim
PO(\bar{F},\alpha\textbf{1},\varphi_1)$. Then from Theorem
\ref{thsdis}, we have $X_{1:n}\leq_{disp}Z_{1:n}$. Again from
Theorem \ref{thsdis2}, we have $Z_{1:n}\leq_{disp}Y_{1:n}$. This
yields $X_{1:n}\leq_{disp}Y_{1:n}$. $\hfill\Box$

The following theorem compare the minimums of two samples, one from
$n$ dependent heterogeneous r.v.'s following the PO model and
another from $n$ dependent homogeneous r.v.'s following the PO
model, in terms of star order.
\begin{thm}\label{thsstar} Suppose $X\sim PO(\bar{F},\boldsymbol\alpha,\varphi)$ and $Y\sim
PO(\bar{F},\alpha\textbf{1},\varphi)$. Then for
$\alpha\geq\frac{1}{n}\sum_{i=1}^n\alpha_{i}$,
$X_{1:n}\leq_{\star}Y_{1:n}$ if $xr(x)$ is decreasing, $\varphi$ is
log-convex, $\frac{\varphi}{\varphi'}$ is concave and $0\leq\alpha
\leq 1$.
\end{thm}
\textbf{Proof:} Using equations (\ref{eqf1s}), (\ref{eqGFs}) and
(\ref{eqgGs}), we have
\begin{eqnarray} \nonumber
&&x^2 \frac{d}{dx}\left(\frac{G_1^{-1}(F_1(x))}{x}\right)\\&=&
\nonumber x \frac{d}{dx}\left(G_1^{-1}(F_1(x))\right)-G_1^{-1}(F_1(x))\\
\nonumber &=& x
\frac{f_{1}(x)}{g_{1}\left(G_1^{-1}(F_1(x))\right)}-G_1^{-1}(F_1(x))\\
\label{ddxGs}&=& \frac{\alpha x r(x) \frac{1}{n}\sum_{i=1}^n
\frac{\varphi\left(\phi\left(\bar{F}_{X_i}(x)\right)\right)}{\varphi'\left(\phi\left(\bar{F}_{X_i}(x)\right)\right)}
\frac{1}{1-\bar{\alpha}_i
\bar{F}(x)}}{r\left(\bar{F}^{-1}(\gamma(x))\right)\left(\alpha
+\bar{\alpha}\varphi\left(\frac{1}{n}\sum_{i=1}^n
\phi\left(\bar{F}_{X_i}(x)\right)\right)\right)
\frac{\varphi\left(\frac{1}{n}\sum_{i=1}^n
\phi\left(\bar{F}_{X_i}(x)\right)\right)}{\varphi'\left(\frac{1}{n}\sum_{i=1}^n
\phi\left(\bar{F}_{X_i}(x)\right)\right)}}-F^{-1}(\gamma(x)).
\end{eqnarray}
In Theorem \ref{thsdis}, for $0\leq\alpha \leq 1$ we have already
proved that
\begin{equation}\label{finvs}
\bar{F}^{-1}(\gamma(x))\geq x.
\end{equation}
Now, if $x r(x)$ is decreasing in $x$, then we have $x r(x)\geq
\bar{F}^{-1}(\gamma(x)) r\left(\bar{F}^{-1}(\gamma(x))\right)$, that
is
\begin{equation}\label{xrxs}
\frac{x r(x)}{r(\bar{F}^{-1}(\gamma(x)))}\geq
\bar{F}^{-1}(\gamma(x)).
\end{equation}
According to the equations (\ref{alpha+}) and (\ref{cheby}) of theorem (\ref{thsdis}), we get
\begin{equation}\label{geq1s} \frac{\frac{\alpha}{n}\sum_{i=1}^n
-\left(
\frac{\varphi\left(\phi\left(\bar{F}_{X_i}(x)\right)\right)}{\varphi'\left(\phi\left(\bar{F}_{X_i}(x)\right)\right)}\right)\frac{1}{1-\bar{\alpha_i}
\bar{F}(x)}}{\left(\alpha
+\bar{\alpha}\varphi\left(\frac{1}{n}\sum_{i=1}^n
\phi\left(\bar{F}_{X_i}(x)\right)\right)\right)
\left(-\frac{\varphi\left(\frac{1}{n}\sum_{i=1}^n
\phi\left(\bar{F}_{X_i}(x)\right)\right)}{\varphi'\left(\frac{1}{n}\sum_{i=1}^n
\phi\left(\bar{F}_{X_i}(x)\right)\right)}\right)}\geq 1.
\end{equation}
Using (\ref{xrxs}) and (\ref{geq1s}), from (\ref{ddxGs}) we get
$$x^2 \frac{d}{dx}\left(\frac{G_1^{-1}(F_1(x))}{x}\right)\geq 0.$$ So,
$\frac{G_1^{-1}(F_1(x))}{x}$ is increasing in $x\geq 0$. Hence
$X_{n:n}\leq_{\star} Y_{n:n}$.$\hfill\Box$

\begin{remark}\normalfont The expression $x r(x)$ is known as the proportional
failure rate (also known as the generalized failure rate)
\citep{rig}. The concerned random variable is said to have the
Decreasing Proportional Failure Rate (DPFR) property if $x r(x)$ is
decreasing in $x$, and in that case domain of the random variable
will be $(0,\infty)$.
\end{remark}

We are interested to know whether as in case of Theorem
\ref{thsstar} we can establish star ordering for $\alpha\geq1$ when
$x r(x)$ is decreasing or increasing. The following counterexample
shows that with these conditions, we cannot establish star ordering
even in case of samples from independent r.v.'s.

\begin{counterexample}
Consider maximums of two samples, one having four independent and
heterogeneous r.v.'s, and another having four independent and
homogeneous r.v.'s. Consider $\alpha_1=0.75$, $\alpha_2=0.95$,
$\alpha_3=23$, $\alpha_4=43$,
$\alpha=(\alpha_1+\alpha_2+\alpha_3+\alpha_4)/4=16.925$, and
$\bar{F}(x)=(1 + \frac{x}{13})^{-0.9}$, so that $x r(x)$ is
increasing. We have
\begin{equation}G_1^{-1}(F_{1}(x))=\bar{F}^{-1}\left(\frac{
\left(\prod_{i=1}^{4}
\bar{F}_{X_i}(x)\right)^{1/4}}{\alpha+\bar{\alpha}\left(\prod_{i=1}^{4}
\bar{F}_{X_i}(x)\right)^{1/4}}\right).\end{equation} We plot
$G_1^{-1}(F_{1}(x))/x$ by substituting $x=t/(1-t)$, so that for
$x\in[0,\infty)$, we have $t\in[0,1)$. From the Figure
\ref{Figsstar}(a), we observe that $G_1^{-1}(F_{1}(x))/x$
is neither increasing nor decreasing.\\
Next we take $\alpha_1=2$, $\alpha_2=33$, $\alpha_3=63$,
$\alpha_4=183$,
$\alpha=(\alpha_1+\alpha_2+\alpha_3+\alpha_4)/4=281/4$, and
$\bar{F}(x)=\frac{1}{x^2}$, $x\in[1,\infty)$ so that $x r(x)$ is
decreasing. We plot $G_1^{-1}(F_{1}(x))/x$ by substituting $x=1/t$,
so that for $x\in[1,\infty)$, we have $t\in[0,1)$. From the Figure
\ref{Figsstar}(b), we observe that $G_1^{-1}(F_{1}(x))/x$ is neither
increasing nor decreasing.

\begin{figure}
\begin{center}
\includegraphics[width=12cm]{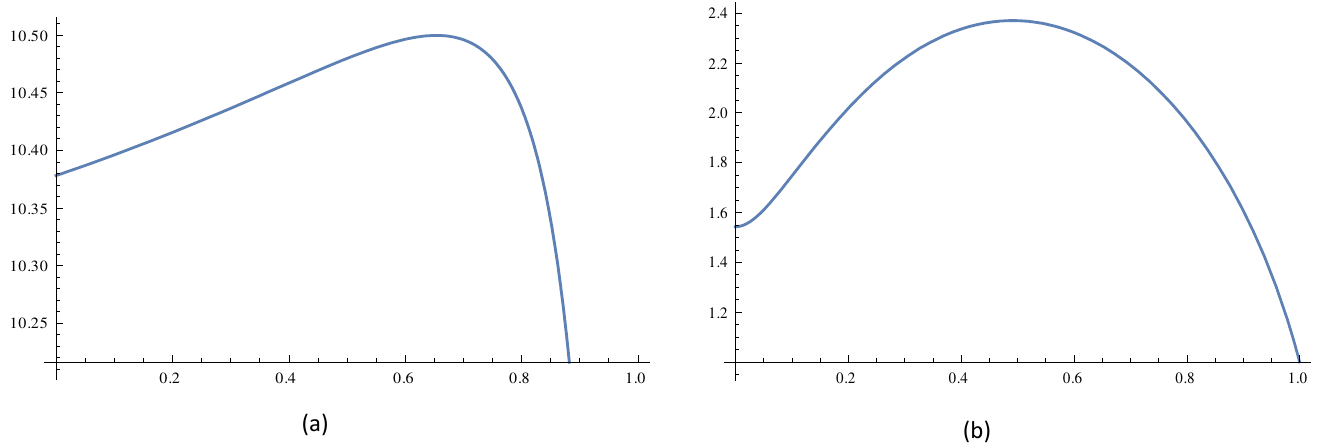}
\caption{Plot of $G_1^{-1}(F_{1}(x))/x$ for (a) $x=t/(1-t)$ when $x
r(x)$ is increasing and (b) $x=1/t$ when $x r(x)$ is decreasing,
$t\in[0,1]$}\label{Figsstar}\end{center}
\end{figure}
\end{counterexample}
The following theorem compares the minimum of two samples, both from
$n$ dependent homogeneous r.v.'s following the PO model and with
different Archimedean copulas. The proof can be done using the
results of proof of Theorem \ref{thsdis2} in the same line as of
Theorem \ref{thsstar}, and hence omitted.
\begin{thm}\label{thsstar2} Suppose $X\sim PO(\bar{F},\alpha\textbf{1},\varphi_1)$ and $Y\sim
PO(\bar{F},\alpha\textbf{1},\varphi_2)$. Then
$X_{1:n}\leq_{\ast}Y_{1:n}$ if $xr(x)$ is decreasing,
$\varphi_2(\phi_2(w)/n)/\varphi_1(\phi_1(w)/n)$ is increasing in
$w$, and $0\leq\alpha \leq 1$.
\end{thm}
The following corollary follows from Theorems \ref{thsstar} and
\ref{thsstar2}.
\begin{coro}Suppose $X\sim PO(\bar{F},\boldsymbol\alpha,\varphi_1)$ and $Y\sim
PO(\bar{F},\alpha\textbf{1},\varphi_2)$. Then for
$\alpha\geq\frac{1}{n}\sum_{i=1}^n\alpha_{i}$,
$X_{1:n}\leq_{\ast}Y_{1:n}$ if $xr(x)$ is decreasing, $\varphi_1$ is
log-convex, $\frac{\varphi_1}{\varphi_1'}$ is concave,
$\varphi_2(\phi_2(w)/n)/\varphi_1(\phi_1(w)/n)$ is increasing in
$w$, and $0\leq\alpha \leq 1$.\end{coro}

\section{The dispersive ordering and the star ordering of maximum of dependent samples following the PO model}

In this section we compare stochastically the maximums of two dependent
samples, one formed from heterogeneous r.v.'s and another from
homogeneous r.v.'s. The distribution function of $X_i$ and $Y_1$
are $F_{X_i}(x)=\frac{ F(x)}{1-\bar{\alpha}_i\bar{F}(x)}$ and
$F_{Y_1}(x)=\frac{ F(x)}{1-\bar{\alpha}\bar{F}(x)}$, respectively,
where $\bar{\alpha}_i=1-\alpha_i$ for $i=1,2,\ldots,n$, and
$\bar{\alpha}=1-\alpha$. The distribution functions of $X_{n:n}$ and
$Y_{n:n}$ are given by
\begin{equation}\label{eqppo}F_{X_{n:n}}(x)=\varphi\left(\sum_{i=1}^n
\phi\left(F_{X_i}(x)\right)\right),\end{equation} and
\begin{equation}\label{eqppo1}F_{Y_{n:n}}(x)=\varphi\left(n
\phi\left(F_{Y_1}(x)\right)\right),\end{equation} where
$\phi(u)=\varphi^{-1}(u)$, $u\in(0,1]$.

The following theorem compare the maximums of two samples, one from $n$ dependent
heterogeneous r.v.'s following the PO model and another from $n$ dependent homogeneous r.v.'s following the PO model, in terms of dispersive order when the baseline distribution function is increasing reversed hazard rate (IRHR). A distribution $F$ is said to be IHRH distribution if the reversed hazard rate $\tilde{r}(\cdot)$ is increasing. If $\tilde{r}(\cdot)$ is decreasing, then $F$ is called decreasing reversed hazard rate (DRHR) distribution.

\begin{thm}\label{thpdis} Suppose $X\sim PO(\bar{F},\boldsymbol\alpha,\varphi)$ and $Y\sim
PO(\bar{F},\alpha\textbf{1},\varphi)$. Then for
$\alpha\geq\frac{1}{n}\sum_{i=1}^n\alpha_{i}$,
$X_{n:n}\geq_{disp}Y_{n:n}$ if the baseline distribution is IRHR,
$\varphi$ is log-concave and $\frac{\varphi}{\varphi'}$ is convex.
\end{thm}
\textbf{Proof:} From equations (\ref{eqppo}) and (\ref{eqppo1}), we
have the distribution functions of $X_{n:n}$ and $Y_{n:n}$
$F_{2}(x)=\varphi\left(\sum_{i=1}^n
\phi\left(F_{X_i}(x)\right)\right)$ and $G_2(x)=\varphi\left(n
\phi\left(F_{Y_1}(x)\right)\right)$, respectively, where
$F_{X_i}(x)=\frac{ F(x)}{\alpha_i+\bar{\alpha}_i F(x)}$ and
$F_{Y_1}(x)=\frac{ F(x)}{\alpha+\bar{\alpha}F(x)}$, $x\in
\mathbb{R}$. The respective pdfs of $X_{n:n} $ and $Y_{n:n}$ are
given by
\begin{equation}\label{eqf2} f_2(x)=\varphi'\left(\sum_{i=1}^n
\phi\left(F_{X_i}(x)\right)\right)\sum_{i=1}^n
\frac{\varphi\left(\phi\left(F_{X_i}(x)\right)\right)}{\varphi'\left(\phi\left(F_{X_i}(x)\right)\right)}
\frac{\alpha_i \tilde{r}(x)}{\alpha_i+\bar{\alpha}_i F(x)},
\end{equation}
\begin{equation*} g_2(x)=n \varphi'\left(n
    \phi\left(F_{Y_1}(x)\right)\right)\cdot \frac{\alpha
        \tilde{r}(x)}{\alpha+\bar{\alpha} F(x)}\cdot
    \frac{\varphi\left(\phi\left(F_{Y_1}(x)\right)\right)}{\varphi'\left(\phi\left(F_{Y_1}(x)\right)\right)},
\end{equation*}
We have
\begin{equation*} G_2^{-1}(x)=F^{-1}\left(\frac{\alpha
\varphi\left(\frac{1}{n}\phi(x)\right)}{1-\bar{\alpha}\varphi\left(\frac{1}{n}\phi(x)\right)}\right),
\end{equation*}
and hence
\begin{equation}\label{eqGF}
G_2^{-1}(F_{2}(x))=F^{-1}\left(\frac{\alpha
\varphi\left(\frac{1}{n}\sum_{i=1}^n
\phi\left(F_{X_i}(x)\right)\right)}{1-\bar{\alpha}\varphi\left(\frac{1}{n}\sum_{i=1}^n
\phi\left(F_{X_i}(x)\right)\right)}\right) =F^{-1}(\beta(x)),
\end{equation} where $\beta(x)=\frac{\alpha
\varphi\left(\frac{1}{n}\sum_{i=1}^n
\phi\left(F_{X_i}(x)\right)\right)}{1-\bar{\alpha}\varphi\left(\frac{1}{n}\sum_{i=1}^n
\phi\left(F_{X_i}(x)\right)\right)}$,
\begin{eqnarray}\nonumber g_{2}(G_2^{-1}(F_{2}(x)))&=&n \varphi'\left(n
\phi\left(\frac{\beta(x)}{\alpha+\bar{\alpha}\beta(x)}\right)\right)\cdot
\frac{\alpha
\tilde{r}\left(\bar{F}^{-1}(\beta(x))\right)}{\alpha+\bar{\alpha}
\beta(x)}\cdot
\frac{\varphi\left(\phi\left(\frac{\beta(x)}{\alpha+\bar{\alpha}\beta(x)}\right)\right)}{\varphi'\left(\phi\left(\frac{
\beta(x)}{\alpha+\bar{\alpha}\beta(x)}\right)\right)}\\
\nonumber &=&n \varphi'\left(\sum_{i=1}^n
\phi\left(F_{X_i}(x)\right)\right) \cdot
\frac{\varphi\left(\frac{1}{n}\sum_{i=1}^n
\phi\left(F_{X_i}(x)\right)\right)}{\varphi'\left(\frac{1}{n}\sum_{i=1}^n
\phi\left(F_{X_i}(x)\right)\right)} \nonumber \\
& & \times  \left(1-\bar{\alpha}\varphi\left(\frac{1}{n}\sum_{i=1}^n
\phi\left(F_{X_i}(x)\right)\right)\right) \cdot \tilde{r}\left(\bar{F}^{-1}(\beta(x))\right).
\label{eqgG}
\end{eqnarray}

Note that $\alpha_i/(\alpha_i+\bar{\alpha}_i F(x))$ is increasing and
concave in $\alpha_i$. It can be seen that
$\phi\left(F_{X_i}(x)\right)$ is increasing and concave in
$\alpha_i$ if $\varphi$ is log-concave. First we take
$\bar{\alpha}\leq 0$. For $\alpha\geq
\frac{1}{n}\sum_{i=1}^n\alpha_{i}$, from the concavity and
increasing property of $\phi\left(F_{X_i}(x)\right)$ with respect to
$\alpha_i$, we have
\begin{eqnarray}
\label{conphip1} \phi\left(F_{Y_1}(x)\right)&\geq&
\frac{1}{n}\sum_{i=1}^n \phi\left(F_{X_i}(x)\right)\\
\label{conphip}
\implies  1- \bar{\alpha} F_{Y_1}(x) &\leq&
1- \bar{\alpha} \varphi\left(\frac{1}{n}\sum_{i=1}^n
\phi\left(F_{X_i}(x)\right)\right)\\
\nonumber
\implies 1-\frac{1}{1- \bar{\alpha}
\varphi\left(\frac{1}{n}\sum_{i=1}^n
\phi\left(F_{X_i}(x)\right)\right)} &\geq&
1-\frac{1}{1- \bar{\alpha} F_{Y_1}(x)}\\
\nonumber
\implies  \frac{\alpha
\varphi\left(\frac{1}{n}\sum_{i=1}^n
\phi\left(F_{X_i}(x)\right)\right)}{1-\bar{\alpha}\varphi\left(\frac{1}{n}\sum_{i=1}^n
\phi\left(F_{X_i}(x)\right)\right)} &\geq&
\frac{\alpha F_{Y_1}(x)}{1-\bar{\alpha}F_{Y_1}(x)}\\
\nonumber
\implies \beta(x)&\geq& F(x).
\end{eqnarray}
Similarly for $\bar{\alpha}\geq 0$, from (\ref{conphip1}) we have
$\beta(x)\geq F(x)$. Thus we have $F^{-1}(\beta(x))\geq x$. Now if
$\tilde{r}(\cdot)$ is increasing then
\begin{equation}\label{eqth11}\tilde{r}(F^{-1}(\beta(x)))\geq \tilde{r}(x).
\end{equation}
Next we take $\bar{\alpha}\geq 0$. As $\varphi(x)$ is decreasing and
convex, we have
\begin{eqnarray}
\nonumber \varphi\left(\frac{1}{n}\sum_{i=1}^n
\phi\left(F_{X_i}(x)\right)\right)&\leq&
\frac{1}{n}\sum_{i=1}^n  \varphi\left(\phi\left(F_{X_i}(x)\right)\right)\\
\label{conphip2}
\implies \bar{\alpha}
\varphi\left(\frac{1}{n}\sum_{i=1}^n
\phi\left(F_{X_i}(x)\right)\right) &\leq& \bar{\alpha} \frac{1}{n}\sum_{i=1}^n \frac{F(x)}{\alpha_i +\bar{\alpha}_i F(x)}\\
\nonumber  &=& \frac{1}{n}\sum_{i=1}^n \bar{\alpha}_i \cdot
\frac{1}{n}\sum_{i=1}^n \frac{F(x)}{\alpha_i +\bar{\alpha}_i F(x)}
\\ \nonumber &\leq& \frac{1}{n}\sum_{i=1}^n \frac{\bar{\alpha}_i F(x)}{\alpha_i +\bar{\alpha}_i F(x)}
\\
\nonumber
\implies 1-\bar{\alpha}
\varphi\left(\frac{1}{n}\sum_{i=1}^n
\phi\left(F_{X_i}(x)\right)\right)&\geq& 1-\frac{1}{n}\sum_{i=1}^n
\frac{\bar{\alpha}_i F(x)}{\alpha_i +\bar{\alpha}_i F(x)}\\\nonumber
&=&\frac{1}{n}\sum_{i=1}^n \frac{\alpha_i }{\alpha_i +\bar{\alpha}_i
F(x)}
\end{eqnarray}
Now for $\bar{\alpha}\leq 0$, from (\ref{conphip}), we have
\begin{eqnarray}\nonumber 1- \bar{\alpha} \varphi\left(\frac{1}{n}\sum_{i=1}^n
\phi\left(F_{X_i}(x)\right)\right) &\geq&
\frac{\alpha}{\alpha+\bar{\alpha}F(x)}\\
\label{eqth13p}&\geq&
\frac{1}{n}\sum_{i=1}^{n}\frac{\alpha_i}{\alpha_i +\bar{\alpha_i}
F(x)},
\end{eqnarray}
where the second inequality follows from the fact that
$\frac{\alpha_i}{\alpha_i+\bar{\alpha}_iF(x)}$ is increasing and
concave in $\alpha_i$.\\ If $\frac{\varphi}{\varphi'}$ is convex,
then we have
\begin{equation}-\frac{\varphi\left(\frac{1}{n}\sum_{i=1}^n
\phi\left(F_{X_i}(x)\right)\right)}{\varphi'\left(\frac{1}{n}\sum_{i=1}^n
\phi\left(F_{X_i}(x)\right)\right)}\geq -\frac{1}{n}\sum_{i=1}^n
\frac{\varphi\left(\phi\left(F_{X_i}(x)\right)\right)}{\varphi'\left(\phi\left(F_{X_i}(x)\right)\right)}.\end{equation}
Thus we have
\begin{equation}\label{alpha+p}\left(1- \bar{\alpha} \varphi\left(\frac{1}{n}\sum_{i=1}^n
\phi\left(F_{X_i}(x)\right)\right)\right)\left(-\frac{\varphi\left(\frac{1}{n}\sum_{i=1}^n
\phi\left(F_{X_i}(x)\right)\right)}{\varphi'\left(\frac{1}{n}\sum_{i=1}^n
\phi\left(F_{X_i}(x)\right)\right)}\right)\geq
\frac{1}{n}\sum_{i=1}^n\left(-
\frac{\varphi\left(\phi\left(F_{X_i}(x)\right)\right)}{\varphi'\left(\phi\left(F_{X_i}(x)\right)\right)}\right)\frac{1}{n}\sum_{i=1}^{n}\frac{\alpha_i}{\alpha_i+\bar{\alpha_i}
F(x)}
\end{equation}
If $\varphi$ is log-concave, then $-\frac{\varphi(x)}{\varphi'(x)}$
is decreasing in $x$, so that
$-\frac{\varphi\left(\phi\left(F_{X_i}(x)\right)\right)}{\varphi'\left(\phi\left(F_{X_i}(x)\right)\right)}$
is decreasing in $\alpha_i$. So by Chebyshev's inequality we have
\begin{equation}\label{chebyp}\frac{1}{n}\sum_{i=1}^n
\left(-\frac{\varphi\left(\phi\left(F_{X_i}(x)\right)\right)}{\varphi'\left(\phi\left(F_{X_i}(x)\right)\right)}\right)\cdot
\frac{1}{n} \sum_{i=1}^n\frac{\alpha_i}{\alpha_i+\bar{\alpha}_i
F(x)} \geq \frac{1}{n}\sum_{i=1}^n
\left(-\frac{\varphi\left(\phi\left(F_{X_i}(x)\right)\right)}{\varphi'\left(\phi\left(F_{X_i}(x)\right)\right)}\right)
\frac{\alpha_i}{\alpha_i+\bar{\alpha}_i F(x)}
\end{equation}
From (\ref{eqth11}), (\ref{alpha+p}), (\ref{chebyp}) and the fact
that the common factor $\varphi'\left(\sum_{i=1}^n
\phi\left(F_{X_i}(x)\right)\right)$ in (\ref{eqf2}) and (\ref{eqgG})
is negative, we have $g_{2}(G_2^{-1}(F_{2}(x)))\geq f_{2}(x)$ for
all $x\in \mathbb{R}$. Hence the theorem follows.$\hfill\Box$

\begin{remark}\normalfont A random variable having support $[0,\infty)$ cannot be IRHR,
however, a distribution function with finite support or the support
of the form $(-\infty,b]$, $0\leq b<\infty$, can be IRHR. For
example, the following distribution functions are IRHR:\\
(i) $F(x)=e^{-(-\lambda x)^\beta}$, $\lambda>0$, $\beta\leq 1$ for
$x\in (-\infty,0]$.\\
(ii) $F(x)=\left(\frac{b q-\mu+x(1-q)}{b-\mu}\right)^{q/(1-q)}$,
$q>1$ for $x\in (-\infty,b]$.\\
It is to be also noted that Archimedean copula with generator
$\psi(x)=\left[1+(2^{-\theta}-1)e^{-x}\right]^{-1/\theta}-1,\text{
where}~ \theta\in(-\infty,0)$ satisfies the condition that $\varphi$
is log-concave and $\frac{\varphi}{\varphi^{\prime}}$ is convex.
\end{remark}
It is of interest to know whether in case of Theorem \ref{thpdis} we
can establish dispersive ordering when baseline distribution is
DRHR. The following counterexample shows that with these conditions,
we cannot establish dispersive ordering even in case of samples from
independent r.v.'s.

\begin{counterexample}
Consider maximums of two samples, one having four independent and
heterogeneous r.v.'s, and another having four independent and
homogeneous r.v.'s with respective distribution functions
$F_{2}(x)=\prod_{i=1}^{r}\left(\frac{F(x)}{1-\bar{\alpha}_i\bar{F}(x)}\right)$
and $G_{2}(x)=\left(\frac{F(x)}{1-\bar{\alpha}\bar{F}(x)}\right)^4$,
where $\alpha_1=0.9$, $\alpha_2=0.95$, $\alpha_3=27$, $\alpha_4=37$,
$\alpha=(\alpha_1+\alpha_2+\alpha_3+\alpha_4)/4=16.4625$, and
$F(x)=1-e^{-(5x)^{0.5}}$, so that the baseline distribution is DRHR.
We obtain \begin{equation*} f_2(x)=\left(\sum_{i=1}^4
\frac{\alpha_i\tilde{r}(x)}{\alpha_i+\bar{\alpha_i}F(x)}\right)
\prod_{i=1}^{4}\left(\frac{F(x)}{\alpha_i+\bar{\alpha_i}F(x)}\right)
\end{equation*} and
\begin{equation*}g_{2}(G_{2}^{-1}(F_{2}(x)))=4 \left(\prod_{i=1}^{4}F_{X_i}(x)\right)
\left(1-\bar{\alpha}\prod_{i=1}^{4}
\left(F_{X_i}(x)\right)^{1/n}\right)\tilde{r}(F^{-1}(\beta(x))),
\end{equation*}
where $\beta(x)=\frac{\alpha \left(\prod_{i=1}^{4}
F_{X_i}(x)\right)^{1/4}}{1-\bar{\alpha} \left(\prod_{i=1}^{4}
F_{X_i}(x)\right)^{1/4}}$.\\
We plot $g_{2}(G_2^{-1}(F_{2}(x)))-f_{2}(x)$ by substituting
$x=t/(1-t)$, so that for $x\in[0,\infty)$, we have $t\in[0,1)$, as shown in Figure \ref{Figpdisp}. It is observed from Figure \ref{Figpdisp} that
$g_{2}(G_2^{-1}(F_{2}(x)))-f_{2}(x)\nleq 0$ and also
$g_{2}(G_2^{-1}(F_{2}(x)))-f_{2}(x)\nleq 0$.
\begin{figure}
\begin{center}
\includegraphics[width=12cm]{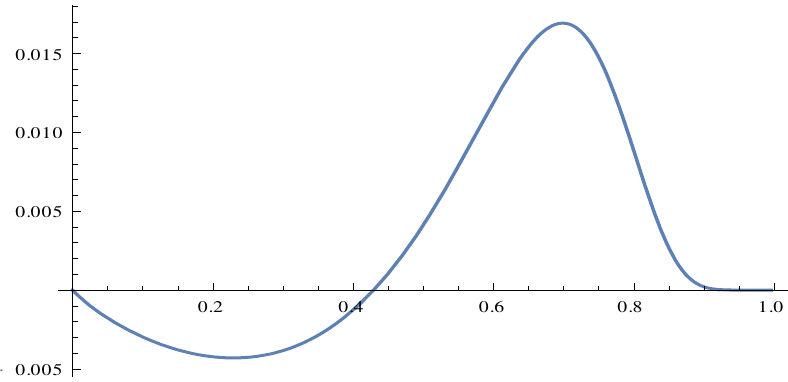}
\caption{Plot of $g_{2}(G_2^{-1}(F_{2}(x)))-f_{2}(x)$ for
$x=t/(1-t)$, $t\in[0,1]$ when baseline distribution is
DRHR.}\label{Figpdisp}\end{center}
\end{figure}
\end{counterexample}

The following theorem compare the maximums of two samples, both from
$n$ dependent homogeneous r.v.'s following the PO model and with
different Archimedean copulas.
\begin{thm}\label{thsdisp2} Suppose $X\sim PO(\bar{F},\alpha\textbf{1},\varphi_1)$ and $Y\sim
PO(\bar{F},\alpha\textbf{1},\varphi_2)$. Then
$X_{n:n}\geq_{disp}Y_{n:n}$ if baseline distribution is IRHR,
$\varphi_1(\phi_1(w)/n)/\varphi_2(\phi_2(w)/n)$ is increasing in
$w$, and $\alpha \geq 1$.
\end{thm}
\textbf{Proof:} The distribution functions of $X_{n:n}$ and
$Y_{n:n}$ are $G_{1}(x)=\varphi_1\left(n
\phi_1\left(F_{X_1}(x)\right)\right)$ and $G_2(x)=\varphi_2\left(n
\phi_2\left(F_{X_1}(x)\right)\right)$, respectively, where
$F_{X_1}(x)=\frac{ F(x)}{\alpha+\bar{\alpha}F(x)}$, $x\in
\mathbb{R}$. We have the pdfs of $X_{n:n}$ and $Y_{n:n}$ are
\begin{equation}\label{eqg1pd} g_1(x)=n \varphi_1'\left(n
\phi_1\left(F_{X_1}(x)\right)\right)\cdot \frac{\alpha
\tilde{r}(x)}{\alpha+\bar{\alpha} F(x)}\cdot
\frac{\varphi_1\left(\phi_1\left(F_{X_1}(x)\right)\right)}{\varphi_1'\left(\phi_1\left(F_{X_1}(x)\right)\right)},
\end{equation}
and
\begin{equation*} g_2(x)=n \varphi_2'\left(n
\phi_2\left(F_{X_1}(x)\right)\right)\cdot \frac{\alpha
\tilde{r}(x)}{\alpha+\bar{\alpha} F(x)}\cdot
\frac{\varphi_2\left(\phi_2\left(F_{X_1}(x)\right)\right)}{\varphi_2'\left(\phi_2\left(F_{X_1}(x)\right)\right)},
\end{equation*}
respectively. We get
\begin{equation}\label{eqGF2}
G_2^{-1}(G_{1}(x))=F^{-1}\left(\frac{\alpha
\varphi_2\left(\frac{1}{n} \phi_2\left(\varphi_1\left(n
\phi_1\left(F_{X_1}(x)\right)\right)\right)\right)}{1-\bar{\alpha}\varphi_2\left(\frac{1}{n}
\phi_2\left(\varphi_1\left(n
\phi_1\left(F_{X_1}(x)\right)\right)\right)\right)}\right)
=F^{-1}(\zeta(x))~(\text{say}),
\end{equation}
\begin{eqnarray}\nonumber g_{2}(G_2^{-1}(G_{1}(x)))&=&n \varphi_2'\left(\phi_2\left(\varphi_1\left(n
\phi_1\left(F_{X_1}(x)\right)\right)\right)\right) \cdot
\frac{\varphi_2\left(\frac{1}{n}\phi_2\left(\varphi_1\left(n
\phi_1\left(F_{X_1}(x)\right)\right)\right)\right)}{\varphi_2'\left(\frac{1}{n}\phi_2\left(\varphi_1\left(n
\phi_1\left(F_{X_1}(x)\right)\right)\right)\right)}\cdot \tilde{r}\left(F^{-1}(\zeta(x))\right) \cdot\\
\label{eqgG2} && \times \left(1-
\bar{\alpha}\varphi_2\left(\frac{1}{n}\phi_2\left(\varphi_1\left(n
\phi_1\left(F_{X_1}(x)\right)\right)\right)\right)\right).\end{eqnarray}
From Lemma 3.9 of \cite{fang1}, for increasing
$\varphi_1(\phi_1(w)/n)/\varphi_2(\phi_2(w)/n)$ we have
$\varphi_2\left(n \phi_2\left(F_{X_1}(x)\right)\right)\leq
\varphi_1\left(n \phi_1\left(F_{X_1}(x)\right)\right)$, which
implies $F_{X_1}(x) \leq \varphi_2\left(\frac{1}{n}
\phi_2\left(\varphi_1\left(n\phi_1\left(F_{X_1}(x)\right)\right)\right)\right)$.
Again from this we get $F(x)\leq \zeta(x)$ which implies
$F^{-1}(\zeta(x))\geq x$. Thus if $\tilde{r}(\cdot)$ is increasing
then
\begin{equation}\label{eqrinvs2}\tilde{r}(F^{-1}(\zeta(x)))\geq \tilde{r}(x).
\end{equation} Also for $\bar{\alpha}\leq 0$, $\varphi_2\left(\frac{1}{n}
\phi_2\left(\varphi_1\left(n\phi_1\left(F_{X_1}(x)\right)\right)\right)\right)\geq
F_{X_1}(x) $ implies
\begin{equation}\label{eqalp+alpb2} 1-\bar{\alpha}\varphi_2\left(\frac{1}{n}
\phi_2\left(\varphi_1\left(n\phi_1\left(F_{X_1}(x)\right)\right)\right)\right)\geq
\frac{\alpha}{\alpha+\bar{\alpha} \bar{F}(x)}.
\end{equation}
Again from Lemma 3.9 of \cite{fang1} by substituting
$w=\varphi_1\left(n \phi_1\left(F_{X_1}(x)\right)\right)$ in
increasing $\varphi_1(\phi_1(w)/n)/\varphi_2(\phi_2(w)/n)$, we get
\begin{equation}\label{remainp2}\frac{\varphi_2'\left(\phi_2\left(\varphi_1\left(n
\phi_1\left(F_{X_1}(x)\right)\right)\right)\right)
\varphi_2\left(\frac{1}{n}\phi_2\left(\varphi_1\left(n
\phi_1\left(F_{X_1}(x)\right)\right)\right)\right)}{\varphi_2'\left(\frac{1}{n}\phi_2\left(\varphi_1\left(n
\phi_1\left(F_{X_1}(x)\right)\right)\right)\right)}\geq
\frac{\varphi_1'\left(n\phi_1\left(F_{X_1}(x)\right)\right)
\varphi_1\left(\phi_1\left(F_{X_1}(x)\right)\right)}{\varphi_1'\left(\phi_1\left(F_{X_1}(x)\right)\right)}.\end{equation}
Now using (\ref{eqrinvs2}), (\ref{eqalp+alpb2}) and
(\ref{remainp2}), we have $g_{2}(G_2^{-1}(G_{1}(x)))\geq g_1(x)$ for
all $x\in \mathbb{R}$. This completes the proof.$\hfill\Box$

\begin{remark}\normalfont It is to be noted that Archimedean copula with generators
\begin{itemize}
    \item [(i)] $\varphi_1(x)=e^{1-(1+x)^(1/\theta_1)}, \theta_1\in (0,\infty)$ and $\varphi_2(x)=e^{\frac{(1-e^x)}{\theta_2}},
    \theta_2\in(0,1)$,
    \item[(ii)] $\varphi_1(x)= e^{\frac{1}{\theta_1\left(1-e^{x}\right)}}$ and $\varphi_2(x)=e^{\frac{1}{\theta_2\left(1-e^{x}\right)}},$ for $0<\theta_2<\theta_1<1,$
\end{itemize}
satisfy the condition that
$\varphi1(\phi1(w)/n)/\varphi2(\phi2(w)/n)$ is increasing in $w$ for
all $n\in Z.$
\end{remark}

The following corollary follows from Theorems \ref{thpdis} and
\ref{thsdisp2}. This corollary compares the minimum of two samples,
one from from $n$ dependent heterogeneous r.v.'s following the PO
model and another from $n$ dependent homogeneous r.v.'s following
the PO model and with different Archimedean copulas.
\begin{coro}Suppose $X\sim PO(\bar{F},\boldsymbol\alpha,\varphi_1)$ and $Y\sim
PO(\bar{F},\alpha\textbf{1},\varphi_2)$. Then for
$\alpha\geq\frac{1}{n}\sum_{i=1}^n\alpha_{i}$,
$X_{n:n}\geq_{disp}Y_{n:n}$ if the baseline distribution is IRHR,
$\varphi_1$ is log-concave, $\frac{\varphi_1}{\varphi'_1}$ is
convex, $\varphi_1(\phi_1(w)/n)/\varphi_2(\phi_2(w)/n)$ is
increasing in $w$, and $\alpha \geq 1$.\end{coro}

The following theorem compares the minimum of two samples, one from $n$ dependent
heterogeneous r.v.'s following the PO model and another from $n$
dependent homogeneous r.v.'s following the PO model, in terms of
star order.
\begin{thm}\label{thpstar} Suppose $X\sim PO(\bar{F},\boldsymbol\alpha,\varphi)$ and $Y\sim
PO(\bar{F},\alpha\textbf{1},\varphi)$. Then for
$\alpha\geq\frac{1}{n}\sum_{i=1}^n\alpha_{i}$,
$X_{n:n}\geq_{\ast}Y_{n:n}$ if $x \tilde{r}(x)$ is increasing in
$x$, $\varphi$ is log-concave, $\frac{\varphi}{\varphi'}$ is convex.
\end{thm}
\textbf{Proof:} Using equations (\ref{eqf2}), (\ref{eqGF}) and
(\ref{eqgG}), we have
\begin{eqnarray} \nonumber
&&x^2 \frac{d}{dx}\left(\frac{G_{2}^{-1}(F_2(x))}{x}\right)\\&=&
\nonumber x \frac{d}{dx}\left(G_{2}^{-1}(F_2(x))\right)-G_{2}^{-1}(F_2(x))\\
\nonumber &=& x
\frac{f_{2}(x)}{g_{2}\left(G_{2}^{-1}(F_2(x))\right)}-G_{2}^{-1}(F_2(x))\\
\label{ddxG}&=& \frac{x \tilde{r}(x) \frac{1}{n}\sum_{i=1}^n
\frac{\varphi\left(\phi\left(F_{X_i}(x)\right)\right)}{\varphi'\left(\phi\left(F_{X_i}(x)\right)\right)}
\frac{\alpha_i}{\alpha_i+\bar{\alpha}_i
F(x)}}{\tilde{r}\left(\bar{F}^{-1}(\beta(x))\right)
\frac{\varphi\left(\frac{1}{n}\sum_{i=1}^n
\phi\left(F_{X_i}(x)\right)\right)}{\varphi'\left(\frac{1}{n}\sum_{i=1}^n
\phi\left(F_{X_i}(x)\right)\right)}\cdot \left(1-
\bar{\alpha}\varphi\left(\frac{1}{n}\sum_{i=1}^n
\phi\left(F_{X_i}(x)\right)\right)\right)}-F^{-1}(\beta(x)).
\end{eqnarray}
In Theorem \ref{thpdis}, we have already proved that
\begin{equation}\label{finv}
F^{-1}(\beta(x))\geq x.
\end{equation}
Now, if $x \tilde{r}(x)$ is increasing in $x$, then we have from
(\ref{finv}), $x \tilde{r}(x)\leq F^{-1}(\beta(x))
\tilde{r}(F^{-1}(\beta(x)))$, that is
\begin{equation}\label{xrx1}
\frac{x \tilde{r}(x)}{\tilde{r}(F^{-1}(\beta(x)))}\leq
F^{-1}(\beta(x)).
\end{equation}
According to the equations (\ref{alpha+p}) and (\ref{chebyp}) of theorem \ref{thsdisp2}, we have
\begin{equation}\label{leq1}\frac{\frac{1}{n}\sum_{i=1}^n
\left(-\frac{\varphi\left(\phi\left(F_{X_i}(x)\right)\right)}{\varphi'\left(\phi\left(F_{X_i}(x)\right)\right)}\right)
\frac{\alpha_i}{\alpha_i+\bar{\alpha}_i F(x)}}{\left(1- \bar{\alpha}
\varphi\left(\frac{1}{n}\sum_{i=1}^n
\phi\left(F_{X_i}(x)\right)\right)\right)\left(-\frac{\varphi\left(\frac{1}{n}\sum_{i=1}^n
\phi\left(F_{X_i}(x)\right)\right)}{\varphi'\left(\frac{1}{n}\sum_{i=1}^n
\phi\left(F_{X_i}(x)\right)\right)}\right)}\leq 1.
\end{equation}
Using (\ref{xrx1}) and (\ref{leq1}), from (\ref{ddxG}) we get $$x^2
\frac{d}{dx}\left(\frac{G_{2}^{-1}(F_2(x))}{x}\right)\leq 0.$$ So,
$\frac{G_{2}^{-1}(F_2(x))}{x}$ is decreasing in $x\geq 0$. Hence
$X_{n:n}\geq_{\star} Y_{n:n}$. $\hfill\Box$\\
The following counterexample shows that we cannot establish star
ordering as in case of Theorem \ref{thpstar} when $x \tilde{r}(x)$
is decreasing or increasing even in case of samples from independent
r.v.'s.

\begin{figure}
\begin{center}
\includegraphics[width=12cm]{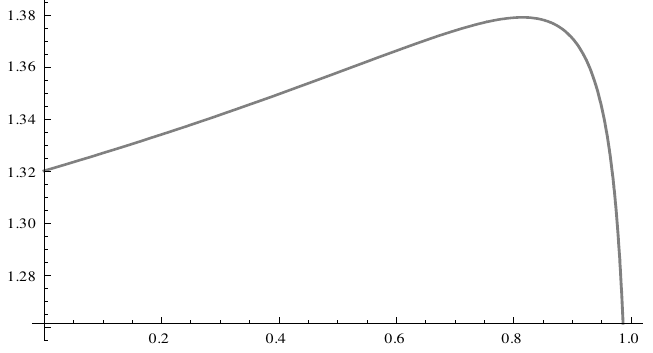}
\caption{Plot of $G_2^{-1}(F_{2}(x))/x$ for $x=t/(1-t)$,
$t\in[0,1]$}\label{Figpstar}\end{center}
\end{figure}

\begin{counterexample}
Consider maximums of two samples, one having four independent and
heterogeneous r.v.'s, and another having four independent and
homogeneous r.v.'s. Consider $\alpha_1=5$, $\alpha_2=15$,
$\alpha_3=25$, $\alpha_4=45$,
$\alpha=(\alpha_1+\alpha_2+\alpha_3+\alpha_4)/4=45/2$, and $F(x)=1 -
(1+x)^{-0.6}$, so that $x \tilde{r}(x)$ is decreasing. We plot
$G_2^{-1}(F_{2}(x))/x$ by substituting $x=t/(1-t)$, so that for
$x\in[0,\infty)$, we have $t\in[0,1)$. We obtain
$$G_{2}^{-1}(F_{2}(x))=F^{-1}\left(\frac{\alpha \left(\prod_{i=1}^{4}
F_{X_i}(x)\right)^{1/4}}{1-\bar{\alpha} \left(\prod_{i=1}^{4}
F_{X_i}(x)\right)^{1/4}}\right).$$ From the Figure \ref{Figpstar},
we observe that $G_2^{-1}(F_{2}(x))/x$ is neither increasing nor
decreasing.
\end{counterexample}

The following theorem compare the minimums of two samples, both from
$n$ dependent homogeneous r.v.'s following the PO model and with
different Archimedean copulas.
\begin{thm}\label{thpstarp2} Suppose $X\sim PO(\bar{F},\alpha\textbf{1},\varphi_1)$ and $Y\sim
PO(\bar{F},\alpha\textbf{1},\varphi_2)$. Then for
$\alpha\geq\frac{1}{n}\sum_{i=1}^n\alpha_{i}$,
$X_{n:n}\geq_{\ast}Y_{n:n}$ if $x \tilde{r}(x)$ is increasing in
$x$, $\varphi_1(\phi_1(w)/n)/\varphi_2(\phi_2(w)/n)$ is increasing
in $w$, and $\alpha\geq 1$.
\end{thm}
\textbf{Proof:} The proof can be done using the
results of proof of Theorem \ref{thsdisp2} in the same line as of
Theorem \ref{thpstar}, and hence omitted.\\

The following corollary follows from Theorems \ref{thpstar} and
\ref{thpstarp2}.
\begin{coro}Suppose $X\sim PO(\bar{F},\boldsymbol\alpha,\varphi_1)$ and $Y\sim
PO(\bar{F},\alpha\textbf{1},\varphi_2)$.Then for
$\alpha\geq\frac{1}{n}\sum_{i=1}^n\alpha_{i}$,
$X_{n:n}\geq_{\ast}Y_{n:n}$ if $x \tilde{r}(x)$ is increasing in
$x$, $\varphi_1$ is log-concave, $\frac{\varphi_1'}{\varphi'}$ is
convex, $\varphi_1(\phi_1(w)/n)/\varphi_2(\phi_2(w)/n)$ is
increasing in $w$, and $\alpha\geq 1$.
\end{coro}

\section{Examples}
Here we demonstrate some of the proposed results numerically. The
first example illustrates the result of Theorem~\ref{thsdis}.

\begin{example}\label{exthsdis}\normalfont
Consider the minimums of two samples, one from three dependent and
heterogeneous r.v.'s, and another from three dependent and
homogeneous r.v.'s, with respective distribution functions
$F_{1}(x)=1-\varphi\left(\sum_{i=1}^3
\phi\left(\frac{\alpha_i\bar{F}(x)}{1-\bar{\alpha}_i\bar{F}(x)}\right)\right)$
and $G_1(x)=1-\varphi\left(3 \phi\left(\frac{ \alpha
\bar{F}(x)}{1-\bar{\alpha}\bar{F}(x)}\right)\right)$, where
$\alpha_1=0.34$, $\alpha_2=0.65$, $\alpha_3=1.23$,
$\alpha=0.88>0.74=(\alpha_1+\alpha_2+\alpha_3)/3$, and
$\bar{F}(x)=e^{-x^{0.3}}$, so that the baseline distribution is DFR.
We take $\varphi(x)=a/\log(x+e^a),~a\in(0,\infty)$ (4.2.19,
\cite{nel}) which satisfies all the conditions of
Theorem~\ref{thsdis}. For this example we take $a=5$. We plot
$g_{1}(G_1^{-1}(F_{1}(x)))-f_{1}(x)$ by substituting $x=t/(1-t)$, so
that for $x\in[0,\infty)$, we have $t\in[0,1)$. The plot is shown in
Figure \ref{sdispex}  and we observe from the plot that
$g_{1}(G_1^{-1}(F_{1}(x)))\leq f_{1}(x)$. Thus
$X_{1:3}\leq_{disp}Y_{1:3}$.

\begin{figure}
\begin{center}
\includegraphics[width=12cm]{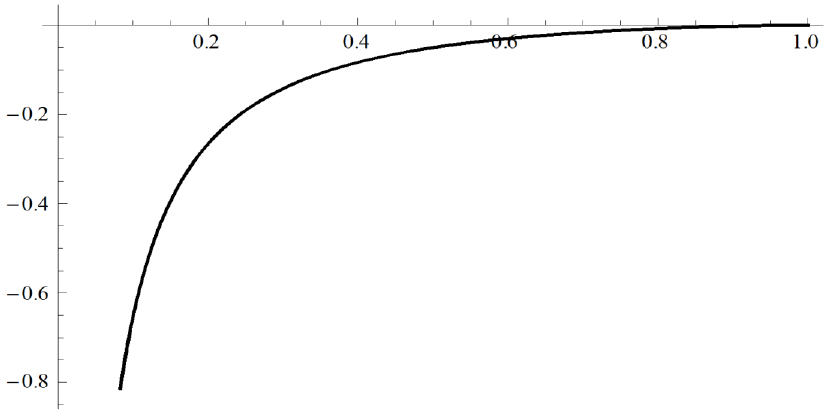}
\caption{Plot of $g_{1}(G_1^{-1}(F_{1}(x)))-f_{1}(x)$ for
$x=t/(1-t)$, $t\in[0,1]$ when baseline distribution is
DFR.}\label{sdispex}\end{center}
\end{figure}
\end{example}

The following example illustrates the result of
Theorem~\ref{thsstar}.

\begin{example} \normalfont
Consider the minimums of two samples, one from four dependent and
heterogeneous r.v.'s, and another from four dependent and
homogeneous r.v.'s. Consider $\alpha_1=0.24$, $\alpha_2=0.45$,
$\alpha_3=0.57$, $\alpha_3=0.57$, $\alpha_4=1.23$,
$\alpha=0.73>(\alpha_1+\alpha_2+\alpha_3+\alpha_4)/4=0.6225$, and
$\bar{F}(x)=1/\sqrt{x}$, $x\in[1,\infty)$ so that $x r(x)$ is
constant. We take $\varphi(x)=a/\log(x+e^a),~a\in(0,\infty)$ which
satisfies all the conditions of Theorem~\ref{thsstar}. For this
example we take $a=7$. We plot $\left(G_1^{-1}(F_{1}(x))/x\right)'$
by substituting $x=1/t$, so that for $x\in[1,\infty)$, we have
$t\in(0,1]$, as shown in Figure \ref{sstarex}. From the figure, we
observe that $G_1^{-1}(F_{1}(x))/x$ is increasing. Thus
$X_{1:4}\leq_{\star}Y_{1:4}$.
\begin{figure}
\begin{center}
\includegraphics[width=12cm]{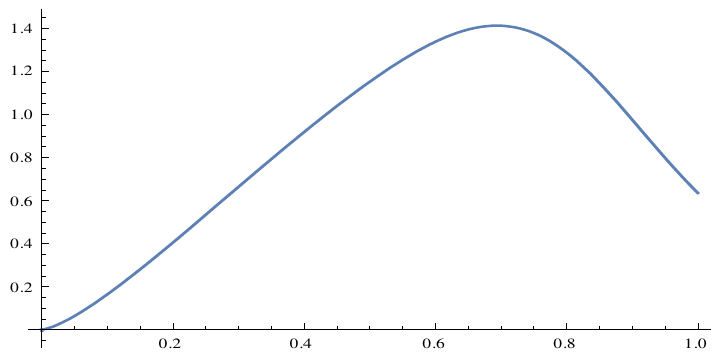}
\caption{Plot of $\left(G_1^{-1}(F_{1}(x))/x\right)'$ for $x=1/t$,
$t\in[0,1] $when $x r(x)$ is decreasing.}\label{sstarex}\end{center}
\end{figure}
\end{example}

The following example illustrates the result of
Theorem~\ref{thpstar}.
\begin{figure}
\begin{center}
\includegraphics[width=12cm]{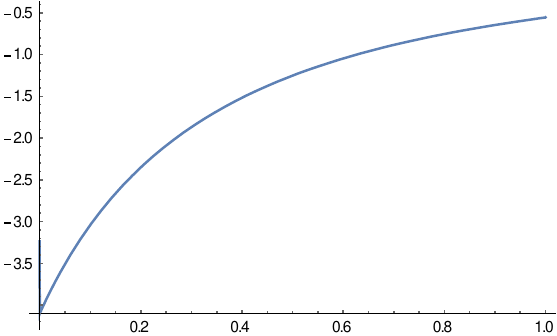}
\caption{Plot of $\left(G_2^{-1}(F_{2}(x))/x\right)'$,
$x\in[0,1]$}\label{pstarex}\end{center}
\end{figure}
\begin{example}\normalfont
Consider the maximums of two samples, one from three dependent and
heterogeneous r.v.'s, and another from three dependent and
homogeneous r.v.'s. Consider $\alpha_1=0.5$, $\alpha_2=0.8$,
$\alpha_3=1.7$, $\alpha=1.6>(\alpha_1+\alpha_2+\alpha_3)/3=1$, and
$F(x)=(e^x-1)/(e-1)$, $x\in[0,1]$ so that $x \tilde{r}(x)$ is
decreasing. We take Archimedean copula with generator
$\psi(x)=\left[1+(2^{-\theta}-1)e^{-x}\right]^{-1/\theta}-1,\text{
with}~ \theta=5$ (4.2.17, \cite{nel}) which satisfies all the
conditions of Theorem~\ref{thpstar}. We plot
$\left(G_2^{-1}(F_{2}(x))/x\right)'$ in Figure \ref{pstarex}. It is
observed from the figure that $G_2^{-1}(F_{2}(x))/x$ is decreasing.
Thus $X_{3:3}\geq_{\star}Y_{3:3}$.
\end{example}

\section{Conclusion}
In this work, we considered the dispersive and the star orders
between both maximums and minimums of samples following the PO model
and coupled with Archimedean copula. The results are illustrated
with numerical examples. Comparing extreme order statistics by means
of some other variability orders or skewness orders like the the
excess wealth order, convex transform order and the Lorenz orders
will be considered in future research.

\bibliographystyle{chicago}
\bibliography{dispersivePO}

\end{document}